\documentclass[10pt]{amsart}
\usepackage{amsmath,amssymb,latexsym,hyperref,url,graphicx}
\usepackage[usenames,dvipsnames]{color}

\graphicspath{{figures/}}

\newtheorem{theorem}{Theorem}
\newtheorem{conjecture}[theorem]{Conjecture}

\newcommand{\eqn}[1]{Equation~(\ref{#1})}
\newcommand{\scn}[1]{Section~\ref{#1}}	

\newcommand{\unabridged}[1]{}
\newcommand{\abridged}[1]{#1}

\newcommand{\depth}{\textnormal{depth}}
\newcommand{\Av}{\textnormal{Av}}
\newcommand{\En}{\textnormal{En}}
\newcommand{\weight}{\textnormal{weight}}
\newcommand{\open}{\textnormal{\texttt{(}}}
\newcommand{\close}{\textnormal{\texttt{)}}}
\newcommand{\gray}[1]{\textcolor{Gray}{#1}}

\newcommand{\vcentergraphics}[1]{\ensuremath{\vcenter{\hbox{\includegraphics{#1}}}}}
\newcommand{\st}[2]{\ensuremath{\vcenter{\hbox{\scalebox{.7}{\includegraphics{binarytree#1-#2}}}}}}

\begin{document}

\title{Pattern avoidance in binary trees}
\author{Eric S. Rowland}
\address{
	Mathematics Department \\
	Tulane University \\
	New Orleans, LA 70118, USA
}
\date{February 8, 2010}

\begin{abstract}
This paper considers the enumeration of trees avoiding a contiguous pattern.  We provide an algorithm for computing the generating function that counts $n$-leaf binary trees avoiding a given binary tree pattern $t$.  Equipped with this counting mechanism, we study the analogue of Wilf equivalence in which two tree patterns are equivalent if the respective $n$-leaf trees that avoid them are equinumerous.  We investigate the equivalence classes combinatorially\unabridged{, finding some relationships to Dyck words avoiding a given subword}.  Toward establishing bijective proofs of tree pattern equivalence, we develop a general method of restructuring trees that conjecturally succeeds to produce an explicit bijection for each pair of equivalent tree patterns.
\end{abstract}

\maketitle
\markboth{Eric Rowland}{Pattern avoidance in binary trees}

\section{Introduction}\label{introduction}

Determining the number of words of length $n$ on a given alphabet that avoid a certain (contiguous) subword is a classical combinatorial problem that can be solved, for example, by the principle of inclusion--exclusion.  An approach to this question using generating functions is provided by the Goulden--Jackson cluster method \cite{goulden-jackson,noonan-zeilberger}, which utilizes only the self-overlaps (or ``autocorrelations'') of the word being considered.  A natural question is ``When do two words have the same avoiding generating function?''  That is, when are the $n$-letter words avoiding (respectively) $w_1$ and $w_2$ equinumerous for all $n$?  The answer is simple:  precisely when their self-overlaps coincide.  For example, the equivalence classes of length-$4$ words on the alphabet $\{0, 1\}$ are as follows.
\[
\begin{array}{cc}
\text{equivalence class} & \text{self-overlap lengths} \\ \hline
\{0001, 0011, 0111, 1000, 1100, 1110\} & \{4\}\\
\{0010, 0100, 0110, 1001, 1011, 1101\} & \{1, 4\} \\
\{0101, 1010\} & \{2, 4\} \\
\{0000, 1111\} & \{1, 2, 3, 4\} \\
\end{array}
\]

In this paper we consider the analogous questions for plane trees.  All trees in the paper are rooted and ordered.  \unabridged{The \emph{depth} of a vertex is the length of the minimal path to that vertex from the root, and $\depth(T)$ is the maximum vertex depth in the tree $T$.}

\begin{figure}
	\includegraphics{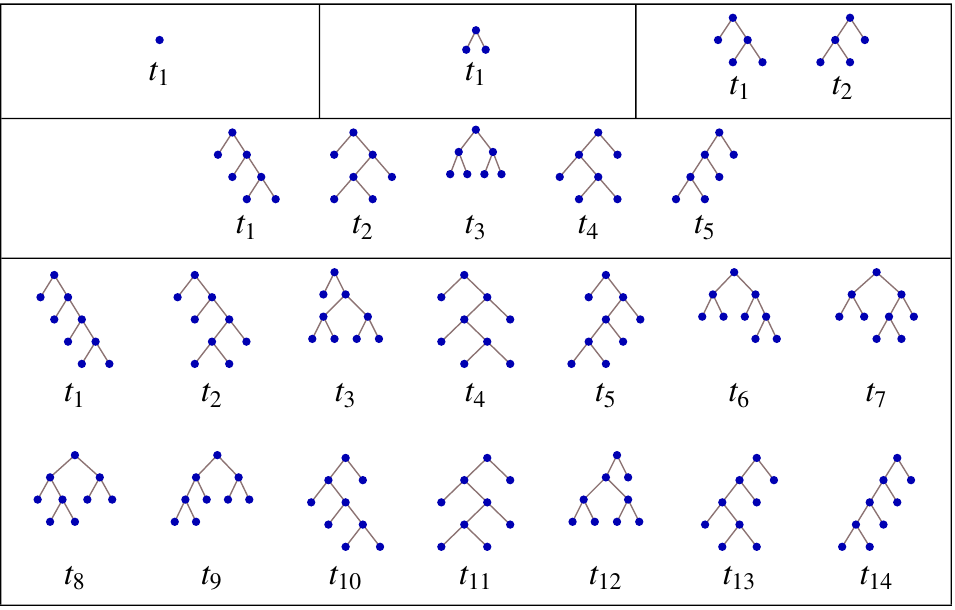}
	\caption{The binary trees with at most $5$ leaves.}\label{binarytrees}
\end{figure}

Our focus will be on \emph{binary trees} --- trees in which each vertex has $0$ or $2$ (ordered) children.  A vertex with $0$ children is a \emph{leaf}, and a vertex with $2$ children is an \emph{internal} vertex.  A binary tree with $n$ leaves has $n-1$ internal vertices, and the number of such trees is the Catalan number $C_{n-1}$.  The first few binary trees are depicted in Figure~\ref{binarytrees}.  We use an indexing for $n$-leaf binary trees that arises from the natural recursive construction of all $n$-leaf binary trees by pairing each $k$-leaf binary tree with each $(n-k)$-leaf binary tree, for all $1 \leq k \leq n-1$.  In practice it will be clear from context which tree we mean by, for example, `$t_1$'.

Conceptually, a binary tree $T$ \emph{avoids} a tree pattern $t$ if there is no instance of $t$ anywhere inside $T$.  Steyaert and Flajolet \cite{sf} were interested in such patterns in vertex-labeled trees.  They were mainly concerned with the asymptotic probability of avoiding a pattern, whereas our focus is on enumeration.  However, they establish in Section~2.2 that the total number of occurrences of an $m$-leaf binary tree pattern $t$ in all $n$-leaf binary trees is
\[
	\binom{2 n - m}{n - m}.
\]
In this sense, all $m$-leaf binary trees are indistinguishable; the results of this paper refine this statement by further distinguishing $m$-leaf tree patterns by the number of $n$-leaf trees containing precisely $k$ copies of each.

We remark that a different notion of tree pattern was later considered by Flajolet, Sipala, and Steyaert \cite{fss}, in which every leaf of the pattern must be matched by a leaf of the tree.  Such patterns are only matched at the bottom of a tree, so they arise naturally in the problem of compactly representing in memory an expression containing repeated subexpressions.  The enumeration of trees avoiding such a pattern is simple, since no two instances of the pattern can overlap:  The number of $n$-leaf binary trees avoiding $t$ depends only on the number of leaves in $t$.  See also Flajolet and Sedgewick \cite[Note III.40]{anacombi}.

The reason for studying patterns in binary trees as opposed to rooted, ordered trees in general is that it is straightforward to determine what it should mean for a binary tree to avoid, for example,
\[
	t_7 = \vcentergraphics{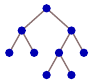},
\]
whereas \emph{a priori} it is ambiguous to say that a general tree avoids
\[
	\vcentergraphics{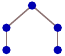}.
\]
Namely, for general trees, `matches a vertex with $i$ children' for $i \geq 1$ could mean either `has exactly $i$ children' or `has at least $i$ children'.  For binary trees, these are the same for $i=2$, so there is no choice to be made.

However, it turns out that the notion of pattern avoidance for binary trees induces a well-defined notion of pattern avoidance for general trees.  This arises via the natural bijection $\beta$ between the set of $n$-leaf binary trees and the set of $n$-vertex trees; using this bijection, one simply translates the problem into the setting of binary trees.

One main theoretical purpose of this paper is to provide an algorithm for computing the generating function that counts binary trees avoiding a certain tree pattern.  This algorithm easily generalizes to count trees containing a prescribed number of occurrences of a certain pattern, and additionally we consider the number of trees containing several patterns each a prescribed number of times.  All of these generating functions are algebraic.  \scn{algorithms} is devoted to these algorithms, which are implemented in \textsc{TreePatterns} \cite{treepatterns}, a \emph{Mathematica} package available from the author's website.

By contrast, another main purpose of this paper is quite concrete, and that is to determine equivalence classes of binary trees.  We say that two tree patterns $s$ and $t$ are \emph{equivalent} if for all $n \geq 1$ the number of $n$-leaf binary trees avoiding $s$ is equal to the number of $n$-leaf binary trees avoiding $t$.  In other words, equivalent trees have the same generating function with respect to avoidance.  This is the analogue of Wilf equivalence in permutation patterns.  Each tree is trivially equivalent to its left--right reflection, but there are other equivalences as well.  The first few classes are presented in \scn{inventory}.  The appendix contains a complete list of equivalence classes of binary trees with at most $6$ leaves, from which we draw examples throughout the paper.  Classes are named with the convention that class~$m.i$ is the $i$th class of $m$-leaf binary trees.

We seek to understand equivalence classes of binary trees combinatorially, and this is the third purpose of the paper.  By analogy with words, one might hope for a simple criterion such as ``$s$ and $t$ are equivalent precisely when the lengths of their self-overlaps coincide''; however, although the set of self-overlap lengths seems to be preserved under equivalence, this statement is not true, for $\{1, 1, 2, 2, 5\}$ corresponds to both classes~6.3 and 6.7.  In lieu of a simple criterion, we look for bijections.  As discussed in \scn{dyckwords}, in a few cases there is a bijection between $n$-leaf binary trees avoiding a certain pattern and Dyck $(n-1)$-words avoiding a certain (contiguous) subword.  In general, when $s$ and $t$ are equivalent tree patterns, we would like to provide a bijection between trees avoiding $s$ and trees avoiding $t$.  Conjecturally, all classes of binary trees can be established bijectively by \emph{top-down} and \emph{bottom-up} replacements; this is the topic of \scn{replacementbijections}.  Nearly all bijections in the paper are implemented in the package \textsc{TreePatterns}.

Aside from mathematical interest, a general study of pattern avoidance in trees has applications to any collection of objects related by a tree structure, such as people in a family tree or species in a phylogenetic tree.  In particular, this paper answers the following question.  Given $n$ related objects (e.g., species) for which the exact relationships aren't known, how likely is it that some prescribed (e.g., evolutionary) relationship exists between some subset of them?  (Unfortunately, it probably will not lead to insight regarding the practical question ``What is the probability of avoiding a mother-in-law?'')  Alternatively, we can think of trees as describing the syntax of sentences in natural language or of fragments of computer code; in this context the paper answers questions about the occurrence and frequency of given phrase substructures.

\abridged{\newpage}	
\section{Definitions}\label{summary}

\subsection{Avoidance}\label{definition}

The more formal way to think of an $n$-vertex tree is as a particular arrangement of $n$ pairs of parentheses, where each vertex is represented by the pair of parentheses containing its children.  For example, the tree
\[
	T = \vcentergraphics{binarytree3-1}
\]
is represented by $\texttt{(()(()()))}$.  This is the word representation of this tree in the alphabet $\{\open, \close\}$.  We do not formally distinguish between the graphical representation of a tree and the word representation, and it is the latter that is useful in manipulating trees algorithmically.  (\emph{Mathematica}'s pattern matching capabilities provide a convenient tool for working with trees represented as nested lists, so this is the convention used by \textsc{TreePatterns}.)

Informally, our concept of containment is as follows.  A binary tree $T$ \emph{contains} $t$ if there is a (contiguous, rooted, ordered) subtree of $T$ that is a copy of $t$.  For example, consider
\[
	t = \vcentergraphics{binarytree3-1}.
\]
None of the trees
\[
	\includegraphics{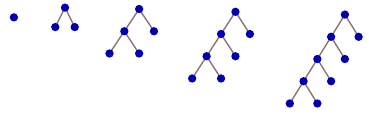}
\]
contains a copy of $t$, while each of the trees
\[
	\includegraphics{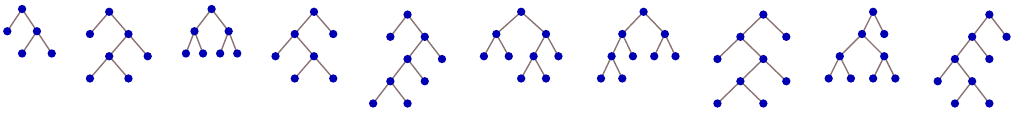}
\]
contains precisely one copy of $t$, each of the trees
\[
	\includegraphics{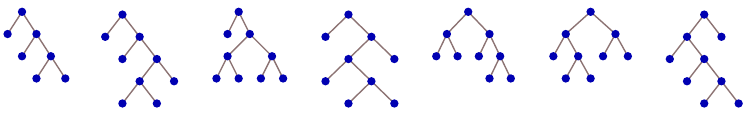}
\]
contains precisely two (possibly overlapping) copies of $t$, and the tree
\[
	\includegraphics{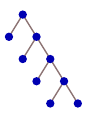}
\]
contains precisely three copies of $t$.  This is a classification of binary trees with at most $5$ leaves according to the number of copies of $t$.

We might formalize this concept with a graph theoretic definition as follows.  Let $t$ be a binary tree.  A \emph{copy} of $t$ in $T$ is a subgraph of $T$ (obtained by removing vertices) that is isomorphic to $t$ (preserving edge directions and the order of children).  Naturally, $T$ \emph{avoids} $t$ if the number of copies of $t$ in $T$ is $0$.

An equivalent but much more useful definition is a language theoretic one, and to provide this we first distinguish a \emph{tree pattern} from a tree.

By `tree pattern', informally we mean a tree whose leaves are ``blanks'' that can be filled (matched) by any tree, not just a single vertex.  More precisely, let $\Sigma = \{\open, \close\}$, and let $L$ be the language on $\Sigma$ containing (the word representation of) every binary tree.  Consider a binary tree $\tau$, and let $t$ be the word on the three symbols $\open, \close, L$ obtained by replacing each leaf $\open\close$ in $\tau$ by $L$.  We call $t$ the \emph{tree pattern} of $\tau$.  This tree pattern naturally generates a language $L_t$ on $\Sigma$, which we obtain by interpreting the word $t$ as a product of the three languages $\open = \{\open\}, \; \close = \{\close\}, \; L$.  Informally, $L_t$ is the set of words that match $t$.  We think of $t$ and $L_t$ interchangeably.  (Note that a tree is a tree pattern matched only by itself.)

For example, let
\[
	\tau = \vcentergraphics{binarytree3-1} = \texttt{(()(()()))};
\]
then the corresponding tree pattern is $t = \open L \open L L \close \close$, and the language $L_t$ consists of all trees of the form $\open T \open U V \close \close$, where $T, U, V$ are binary trees.

Let $\Sigma^*$ denote the set of all finite words on $\Sigma$.  The language $\Sigma^* L_t \Sigma^* \cap L$ is the set of all binary trees whose word has a subword in $L_t$.  Therefore we say that a binary tree $T$ \emph{contains} the tree pattern $t$ if $T$ is in the language $\Sigma^* L_t \Sigma^* \cap L$.  We can think of this language as a multiset, where a given tree $T$ occurs with multiplicity equal to the number of ways that it matches $\Sigma^* L_t \Sigma^*$.  Then the \emph{number of copies} of $t$ in $T$ is the multiplicity of $T$ in $\Sigma^* L_t \Sigma^* \cap L$.

Continuing the example from above, the tree
\[
	T = \vcentergraphics{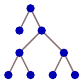} = \texttt{(()((()())(()())))}
\]
contains $2$ copies of $t$ since it matches $\Sigma^* L_t \Sigma^*$ in $2$ ways:  $\open T \open U V \close \close$ with $T = \open \close$ and $U = V = \texttt{(()())}$, and $\gray{\open \open \close} \open T \open U V \close \close \gray{\close}$ with $T = \texttt{(()())}$ and $U = V = \open \close$.

Our notation distinguishes tree patterns from trees:  Tree patterns are represented by lowercase variables, and trees are represented by uppercase variables.  To be absolutely precise, we would graphically distinguish between terminal leaves $\open\close$ of a tree and blank leaves $L$ of a tree pattern, but this gets in the way of speaking about them as the same objects, which is rather convenient.

In Sections~\ref{algorithms} and \ref{replacementbijections} we will be interested in taking the intersection $p \cap q$ of tree patterns $p$ and $q$ (by which we mean the intersection of the corresponding languages $L_p$ and $L_q$).  The intersection of two or more explicit tree patterns can be computed recursively:  $p \cap L = p$, and $\open p_l p_r \close \cap \open q_l q_r \close = \open (p_l \cap q_l) (p_r \cap q_r) \close$.

\subsection{Generating functions}\label{generatingfunctions}

Our primary goal is to determine the number $a_n$ of binary trees with $n$ vertices that avoid a given binary tree pattern $t$, and more generally to determine the number $a_{n,k}$ of binary trees with $n$ vertices and precisely $k$ copies of $t$.  Thus we consider two generating functions associated with $t$: the \emph{avoiding generating function}
\[
	\Av_t(x) = \sum_{\text{$T$ avoids $t$}} x^{\text{number of vertices in $T$}} = \sum_{n=0}^\infty a_n x^n
\]
and the \emph{enumerating generating function}
\begin{align}
	\En_{L, t}(x,y)
	&= \sum_{T \in L} x^{\text{number of vertices in $T$}} y^{\text{number of copies of $t$ in $T$}} \notag \\
	&= \sum_{\vphantom{k} n=0}^\infty \sum_{k=0}^\infty a_{n,k} x^n y^k. \notag
\end{align}
The avoiding generating function is the special case $\Av_t(x) = \En_{L, t}(x,0)$.

\begin{theorem}\label{enumerating}
$\En_{L, t}(x,y)$ is algebraic.
\end{theorem}

The proof is constructive, so it enables us to compute $\En_{L, t}(x)$, and in particular $\Av_t(x)$, for explicit tree patterns.  We postpone the proof until \scn{algorithm-enumeratingsingle} to address a natural question that arises:  Which trees have the same generating function?  That is, for which pairs of binary tree patterns $s$ and $t$ are the $n$-leaf trees avoiding (or containing $k$ copies of) these patterns equinumerous?

We say that $s$ and $t$ are \emph{avoiding-equivalent} if $\Av_s(x) = \Av_t(x)$.  We say they are \emph{enumerating-equivalent} if the seemingly stronger condition $\En_{L, s}(x,y) = \En_{L, t}(x,y)$ holds.  We can compute these equivalence classes explicitly by computing $\Av_t(x)$ and $\En_{L, t}(x,y)$ for, say, all $m$-leaf binary tree patterns $t$.  In doing this for binary trees with up to $7$ leaves, one comes to suspect that these conditions are in fact the same.

\begin{conjecture}
If $s$ and $t$ are avoiding-equivalent, then they are also enumerating-equivalent.
\end{conjecture}

In light of this experimental result, we focus attention in the remainder of the paper on classes of avoiding-equivalence, since conjecturally they are the same as classes of enumerating-equivalence.

\section{Initial inventory and some special bijections}\label{inventory}

In this section we undertake an analysis of small patterns.  We determine $\Av_t(x)$ for binary tree patterns with at most $4$ leaves using methods specific to each.  This allows us to establish the equivalence classes in this range.

\subsection{$1$-leaf trees}\label{1leafclasses}

There is only one binary tree pattern with a single leaf, namely
\[
	t = \vcentergraphics{binarytree1-1} = L.
\]
Every binary tree contains at least one vertex, so $\Av_t(x) = 0$.  The number of binary trees with $2n-1$ vertices is $C_{n-1}$, so
\[
	\En_L(x) = x + x^3 + 2 x^5 + 5 x^7 + 14 x^9 + 42 x^{11} + \cdots = \sum_{n=1}^\infty C_{n-1} x^{2n-1}.
\]

\subsection{$2$-leaf trees}\label{2leafclasses}

There is also only one binary tree pattern with precisely $2$ leaves:
\[
	t = \vcentergraphics{binarytree2-1} = \open L L \close.
\]
However, $t$ is a fairly fundamental structure in binary trees; the only tree avoiding it is the $1$-vertex tree $\open\close$.  Thus $\Av_t(x) = x$, and
\[
	\En_{L, t}(x,y) = \sum_{n=1}^\infty C_{n-1} x^{2n-1} y^{n-1} = \frac{1 - \sqrt{1 - 4 x^2 y}}{2 x y}.
\]

\subsection{$3$-leaf trees}\label{3leafclasses}

There are $C_2 = 2$ binary trees with $3$ leaves, and they are equivalent by left--right reflection:
\[
	\vcentergraphics{binarytree3-1} \quad \text{and} \quad \vcentergraphics{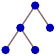}.
\]
There is only one binary tree with $n$ leaves avoiding
\[
	\vcentergraphics{binarytree3-2} = \open \open L L \close L \close,
\]
namely the ``right comb'' $\texttt{(()(()(()(()}\cdots\texttt{))))}$.  Therefore for these trees
\[
	\Av_t(x) = x + x^3 + x^5 + x^7 + x^9 + x^{11} + \cdots = \frac{x}{1 - x^2}.
\]

\subsection{$4$-leaf trees}\label{4leafclasses}

Among $4$-leaf binary trees we find more interesting behavior.  There are $C_3 = 5$ such trees, pictured as follows.
\[
	\includegraphics{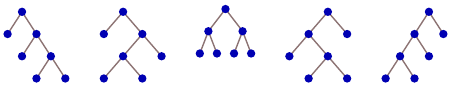}
\]
They comprise $2$ equivalence classes.

\subsubsection*{Class~4.1}\label{binaryclass4-1}

The first equivalence class consists of the trees
\[
	t_1 = \vcentergraphics{binarytree4-1} \quad \text{and} \quad
	t_5 = \vcentergraphics{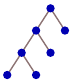}.
\]
The avoiding generating function $\Av_t(x)$ for each of these trees satisfies
\[
	x^3 f^2 + (x^2 - 1) f + x = 0
\]
because the number of $n$-leaf binary trees avoiding $t_1$ is the Motzkin number $M_{n-1}$:
\[
	\Av_t(x) = x + x^3 + 2 x^5 + 4 x^7 + 9 x^9 + 21 x^{11} + \cdots = \sum_{n=1}^\infty M_{n-1} x^{2n-1}.
\]

This fact is presented by Donaghey and Shapiro \cite{donaghey-shapiro} as their final example of objects counted by the Motzkin numbers.  They provide a bijective proof which we reformulate here.  Specifically, there is a natural bijection between the set of $n$-leaf binary trees avoiding $t_1$ and the set of Motzkin paths of length $n-1$ --- paths from $(0,0)$ to $(n-1,0)$ composed of steps $\langle 1, -1 \rangle$, $\langle 1, 0 \rangle$, $\langle 1, 1 \rangle$ that do not go below the $x$-axis.  We represent a Motzkin path as a word on $\{-1, 0, 1\}$ encoding the sequence of steps under $\langle 1, \Delta y \rangle \mapsto \Delta y$.  \unabridged{The bijection is as follows.}

\abridged{Let $\beta$ be the usual bijection from $n$-leaf binary trees to $n$-vertex general trees that operates by contracting every rightward edge.}\unabridged{Recall the map $\beta$ from \scn{harary-prints-tutte}.}  To obtain the Motzkin path associated with a binary tree $T$ avoiding $t_1$:
\begin{enumerate}
	\item Let $T' = \beta(T)$.  No vertex in $T'$ has more than $2$ children, since
	\[
		\beta(t_1) = \vcentergraphics{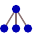}
	\]
	and $T$ avoids $t_1$.
	\item Create a word $w$ on $\{-1, 0, 1\}$ by traversing $T'$ in depth-first order (i.e., for each subtree visit first the root vertex and then its children trees in order); for each vertex, record $1$ less than the number of children of that vertex.
	\item Delete the last letter of $w$ (which is $-1$).
\end{enumerate}
The resulting word contains the same number of $-1$s and $1$s, and every prefix contains at least as many $1$s as $-1$s, so it is a Motzkin path.  The steps are easily reversed to provide the inverse map from Motzkin paths to binary trees avoiding $t_1$.  (For the larger context of this bijection, see Stanley's presentation leading up to Theorem 5.3.10 \cite{Stanley 2}.)

\subsubsection*{Class~4.2}\label{binaryclass4-2}

The second equivalence class consists of the three trees
\[
	t_2 = \vcentergraphics{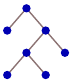}, \quad
	t_3 = \vcentergraphics{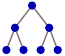}, \quad \text{and} \quad
	t_4 = \vcentergraphics{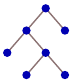}
\]
and provides the smallest example of nontrivial equivalence.  Symmetry gives $\Av_{t_2}(x) = \Av_{t_4}(x)$.  To establish $\Av_{t_2}(x) = \Av_{t_3}(x)$, for each of these trees $t$ we give a bijection between $n$-leaf binary trees avoiding $t$ and binary words of length $n-2$.  By composing these two maps we obtain a bijection between trees avoiding $t_2$ and trees avoiding $t_3$.

First consider
\[
	t_3 = \vcentergraphics{binarytree4-3}.
\]
If $T$ avoids $t_3$, then no vertex of $T$ has four grandchildren; that is, at most one of a vertex's children has children of its own.  This implies that at each generation at most one vertex has children.  Since there are two vertices at each generation after the first, the number of such $n$-leaf trees is $2^{n-2}$ for $n \ge 2$:
\[
	\Av_{t_3}(x) = x + x^3 + 2 x^5 + 4 x^7 + 8 x^9 + 16 x^{11} + \cdots = x + \sum_{n=2}^\infty 2^{n-2} x^{2n-1} = \frac{x (1 - x^2)}{1 - 2 x^2}.
\]
Form a word $w \in \{0,1\}^{n-2}$ corresponding to $T$ by letting the $i$th letter be $0$ or $1$ depending on which vertex (left or right) on level $i+1$ has children.

Now consider
\[
	t_2 = \vcentergraphics{binarytree4-2}.
\]
A ``typical'' binary tree avoiding $t_2$ looks like
\[
	\vcentergraphics{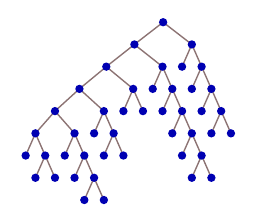}
\]
and is determined by the length of its spine and the length of each arm.  Starting from the root, traverse the internal vertices of a tree $T$ avoiding $t_2$ according to the following rule.  Always move to the right child of a vertex when the right child is an internal vertex, and if the right child is a leaf then move to the highest unvisited internal spine vertex.  By recording $0$ and $1$ for left and right movements in this traversal, a word $w$ on $\{0,1\}$ is produced that encodes $T$ uniquely.  We have $|w| = n-2$ since we obtain one symbol from each internal vertex except the root.  Since every word $w$ corresponds to an $n$-leaf binary tree avoiding $t_2$, there are $2^{n-2}$ such trees.

More formally, let $\omega$ be a map from binary trees to binary words defined by $\omega(\open T_l T_r \close) = \kappa_1(T_r) \kappa_0(T_l)$, where
\[
	\kappa_i(T) =
	\begin{cases}
		\epsilon		&\text{if $T = \open\close$}\\
		i \, \omega(T)		&\text{otherwise.}
	\end{cases}
\]
Then the word corresponding to $T$ is $w = \omega(T)$.

For the inverse map $\omega^{-1}$, begin with the word $\open l r \close$.  Then read $w$ left to right.  When the symbol $1$ is read, replace the existing $r$ by $\open\open\close r \close$; when $0$ is read, replace the existing $r$ by $\open\close$ and the existing $l$ by $\open l r \close$.  After the entire word is read, replace the remaining $l$ and $r$ with $\open\close$.  One verifies that $T$ has $n$ leaves.  The tree $T$ avoids $t_2$ because the left child of an $r$ vertex never has children of its own.

\subsection{Bijections to Dyck words}\label{dyckwords}

\abridged{
We mention that in some cases the set of trees avoiding a pattern is in bijection to the set of Dyck words avoiding a certain subword.  For example, trees avoiding a pattern in class~4.1 are in bijection to Dyck words avoiding the subword $000$.  Recall $\beta$ from \scn{4leafclasses}.  For $t_5$ we have
\[
	\beta\left(\vcentergraphics{binarytree4-5}\right) = \vcentergraphics{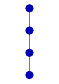},
\]
whose corresponding Dyck word is $000111$, and $\beta$ has the feature that $T$ contains $t_5$ if and only if the Dyck word corresponding to $\beta(T)$ contains $000$.

In general there is a bijection between $n$-leaf binary trees avoiding $t$ and $(n-1)$-Dyck words avoiding $w$ whenever $w$ is a characteristic feature of $\beta(t)$, that is, some feature of the tree that is preserved locally by $\beta$.

For example, for $t_2$ in class~4.2 we observe that
\[
	\beta\left(\vcentergraphics{binarytree4-2}\right) = \vcentergraphics{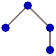} = 010011,
\]
and containing $100$ as a subword is necessary and sufficient for the corresponding tree to contain $t_2$.  Thus binary trees avoiding a pattern in class~4.2 are in bijection to Dyck words avoiding $100$.  Other bijections can be found similarly:  Classes~5.2, 5.3, 6.3, and 6.6 correspond to the words $1100$, $1000$, $11000$, and $10000$.

Notes on sequences counting Dyck words avoiding a subword have been contributed by David Callan and Emeric Deutsch to the Encyclopedia of Integer Sequences \cite{eis}.  The subject appears to have begun with Deutsch \cite[Section~6.17]{deutsch}, who enumerated Dyck words according to the number of occurrences of $100$.  Sapounakis, Tasoulas, and Tsikouras \cite{stt} have considered additional subwords.  Via the bijections just described, their results provide additional derivations of the generating functions $\Av_t(x)$.
}

\section{Algorithms}\label{algorithms}

In this section we provide algorithms for computing algebraic equations satisfied by $\Av_t(x)$, $\En_{L, t}(x,y)$, and the more general $\En_{L, p_1, \dots, p_k}(x_L, x_{p_1}, \dots, x_{p_k})$ defined in \scn{algorithm-multiple}.  Computing $\Av_t(x)$ or $\En_{L, t}(x,y)$ for all $m$-leaf binary tree patters $t$ allows one to automatically determine the equivalence classes given in the appendix.

We draw upon the notation introduced in \scn{definition}.  In particular, the intersection $p \cap p'$ of two tree patterns plays a central role.  Recall that $L_p$ is the set of trees matching $p$ at the top level.

\abridged{The \emph{depth} of a vertex in a tree is the length of the minimal path to that vertex from the root, and $\depth(T)$ is the maximum vertex depth in the tree $T$.}

\subsection{Avoiding a single tree}\label{algorithm-avoidingsingle}

Fix a binary tree pattern $t$ we wish to avoid.  For a given tree pattern $p$, we will make use of the generating function
\[
	\weight(p) = \weight(L_p) := \sum_{T \in L_p} \weight(T),
\]
where
\[
	\weight(T) =
	\begin{cases}
		x^\text{number of vertices in $T$}	&\text{if $T$ avoids $t$}\\
		0					&\text{if $T$ contains $t$.}
	\end{cases}
\]

The case $t = L$ was covered in \scn{1leafclasses}, so we assume $t \neq L$.  Then $t = \open t_l t_r \close$ for some tree patterns $t_l$ and $t_r$.  Since $\open T_l T_r \close$ matches $t$ precisely when $T_l$ matches $t_l$ and $T_r$ matches $t_r$, we have
\begin{multline}\label{logicrule}
	\weight(\open p_l p_r \close) = \\
	x \cdot \big( \weight(p_l) \cdot \weight(p_r) - \weight(p_l \cap t_l) \cdot \weight(p_r \cap t_r) \big).
\end{multline}
The coefficient $x$ is the weight of the root vertex of $\open p_l p_r \close$ that we destroy in separating this pattern into its two subpatterns.

We now construct a polynomial (with coefficients that are polynomials in $x$) that is satisfied by $\Av_t(x) = \weight(L)$, the weight of the language of binary trees.  The algorithm is as follows.

Begin with the equation
\[
	\weight(L) = \weight(\open\close) + \weight(\open L L \close).
\]
The variable $\weight(\open L L \close)$ is ``new''; we haven't yet written it in terms of other variables.  So use \eqn{logicrule} to rewrite $\weight(\open L L \close)$.  For each expression $\weight(p \cap p')$ that is introduced, we compute the intersection $p \cap p'$.  This allows us to write $\weight(p \cap p')$ as $\weight(q)$ for some pattern $q$ that is simply a word on $\{\open, \close, L\}$ (i.e., does not contain the $\cap$ operator).

For each new variable $\weight(q)$, we obtain a new equation by making it the left side of \eqn{logicrule}, and then as before we eliminate $\cap$ by explicitly computing intersections.

We continue in this manner until there are no new variables produced.  This must happen because $\depth(p \cap p') = \max(\depth(p), \depth(p'))$, so since there are only finitely many trees that are shallower than $t$, there are only finitely many variables in this system of polynomial equations.

Finally, we compute a Gr\"obner basis for the system in which all variables except $\weight(\open\close) = x$ and $\weight(L) = \Av_t(x)$ are eliminated.  This gives a single polynomial equation in these variables, establishing that $\Av_t(x)$ is algebraic.

Let us work out an example.  We use the graphical representation of tree patterns with the understanding that the leaves are blanks.  Consider the tree pattern
\[
	t = \vcentergraphics{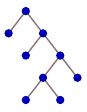} = \open L \open L \open \open L L \close L \close \close \close
\]
from class~5.2.  The first equation is
\[
	\weight(\st{1}{1}) = x + \weight(\st{2}{1}).
\]
We have $t_l = \st{1}{1}$ and $t_r = \st{4}{2}$, so \eqn{logicrule} gives
\begin{align}
	\weight(\st{2}{1})
	&= x \cdot \big( \weight(\st{1}{1}) \cdot \weight(\st{1}{1}) - \weight(\st{1}{1} \cap \st{1}{1}) \cdot \weight(\st{1}{1} \cap \st{4}{2}) \big) \notag \\
	&= x \cdot \big( \weight(\st{1}{1})^2 - \weight(\st{1}{1}) \cdot \weight(\st{4}{2}) \big) \notag
\end{align}
since $L \cap p = p$ for any tree pattern $p$.  The variable $\weight(\st{4}{2}) = \weight(t_r)$ is new, so we put it into \eqn{logicrule}:
\begin{align}
	\weight(\st{4}{2})
	&= x \cdot \big( \weight(\st{1}{1}) \cdot \weight(\st{3}{2}) - \weight(\st{1}{1} \cap \st{1}{1}) \cdot \weight(\st{3}{2} \cap \st{4}{2}) \big) \notag \\
	&= x \cdot \big( \weight(\st{1}{1}) \cdot \weight(\st{3}{2}) - \weight(\st{1}{1}) \cdot \weight(\st{5}{7}) \big). \notag
\end{align}
There are two new variables:
\begin{align}
	\weight(\st{3}{2})
	&= x \cdot \big( \weight(\st{2}{1}) \cdot \weight(\st{1}{1}) - \weight(\st{2}{1} \cap \st{1}{1}) \cdot \weight(\st{1}{1} \cap \st{4}{2}) \big) \notag \\
	&= x \cdot \big( \weight(\st{2}{1}) \cdot \weight(\st{1}{1}) - \weight(\st{2}{1}) \cdot \weight(\st{4}{2}) \big); \notag \\
	\weight(\st{5}{7})
	&= x \cdot \big( \weight(\st{2}{1}) \cdot \weight(\st{3}{2}) - \weight(\st{2}{1} \cap \st{1}{1}) \cdot \weight(\st{3}{2} \cap \st{4}{2}) \big) \notag \\
	&= x \cdot \big( \weight(\st{2}{1}) \cdot \weight(\st{3}{2}) - \weight(\st{2}{1}) \cdot \weight(\st{5}{7}) \big). \notag
\end{align}
We have no new variables, so we eliminate the four auxiliary variables
\[
	\weight(\st{2}{1}), \weight(\st{3}{2}), \weight(\st{4}{2}), \weight(\st{5}{7})
\]
from this system of five equations to obtain
\[
	x^3 \, \weight(\st{1}{1})^2 - (x^2 - 1)^2 \, \weight(\st{1}{1}) - x \, (x^2 - 1) = 0.
\]

\subsection{Enumerating with respect to a single tree}\label{algorithm-enumeratingsingle}

To prove Theorem~\ref{enumerating}, we make a few modifications in order to compute $\En_{L, t}(x,y)$ instead of $\Av_t(x)$.  Again
\[
	\weight(p) := \sum_{T \in L_p} \weight(T),
\]
but now $\weight(T) = x^\text{number of vertices in $T$} y^\text{number of copies of $t$ in $T$}$ for all $T$.  We modify \eqn{logicrule} to become
\begin{multline}\label{enlogicrule}
	\weight(\open p_l p_r \close) = \\
	x \cdot \big(\weight(p_l) \cdot \weight(p_r) + (y-1) \cdot \weight(p_l \cap t_l) \cdot \weight(p_r \cap t_r) \big)
\end{multline}
since in addition to accounting for the trees that avoid $t$ we also account for those that match $t$, in which case $y$ is contributed.

The rest of the algorithm carries over unchanged, and we obtain a polynomial equation in $x$, $y$, and $\En_{L, t}(x,y) = \weight(L)$.

\subsection{Enumerating with respect to multiple trees}\label{algorithm-multiple}

A more general question is the following.  Given several binary tree patterns $p_1, \dots, p_k$, what is the number $a_{n_0,n_1,\dots,n_k}$ of binary trees containing precisely $n_0$ vertices, $n_1$ copies of $p_1$, \dots, $n_k$ copies of $p_k$?  We consider the enumerating generating function
\begin{align}
	\En_{L, p_1, \dots, p_k}(x_L, x_{p_1}, \dots, x_{p_k})
	&= \sum_{T \in L} x_L^{\alpha_0} x_{p_1}^{\alpha_1} \cdots x_{p_k}^{\alpha_k} \notag \\
	&= \sum_{n_0=0}^\infty \sum_{n_1=0}^\infty \cdots \sum_{n_k=0}^\infty a_{n_0,n_1,\dots,n_k} x_L^{n_0} x_{p_1}^{n_1} \cdots x_{p_k}^{n_k}, \notag
\end{align}
where $p_0 = L$ and $\alpha_i$ is the number of copies of $p_i$ in $T$.  (We need not assume that the $p_i$ are distinct.)  This generating function can be used to obtain information about how correlated a family of tree patterns is.  We have the following generalization of Theorem~\ref{enumerating}.

\begin{theorem}\label{enumerating-multiple}
$\En_{L, p_1, \dots, p_k}(x_L, x_{p_1}, \dots, x_{p_k})$ is algebraic.
\end{theorem}

Keeping track of multiple tree patterns $p_1, \dots, p_k$ is not much more complicated than handling a single pattern, and the algorithm for doing so has the same outline.  Let
\[
	\weight(p) := \sum_{T \in L_p} \weight(T)
\]
with
\[
	\weight(T) = x_L^{\alpha_0} x_{p_1}^{\alpha_1} \cdots x_{p_k}^{\alpha_k},
\]
where again $\alpha_i$ is the number of copies of $p_i$ in $T$.  Let $d = \max_{1 \leq i \leq k} \depth(p_i)$.  First we describe what to do with each new variable $\weight(q)$ that arises.  The approach used is different than that for one tree pattern; in particular, we do not make use of intersections.  Consequently, it is less efficient.

Let $l$ be the number of leaves in $q$.  If $T$ is a tree matching $q$, then for each leaf $L$ of $q$ there are two possibilities:  Either $L$ is matched by a terminal vertex $\open\close$ in $T$, or $L$ is matched by a tree matching $\open L L \close$.  For each leaf we make this choice independently, thus partitioning the language $L_q$ into $2^l$ disjoint sets represented by $2^l$ tree patterns that are disjoint in the sense that each tree matching $q$ matches precisely one of these patterns.  For example, partitioning the pattern $\open L L \close$ into $2^2$ patterns gives
\begin{multline}
	\weight(\open L L \close) = \weight(\open \open\close \open\close \close) + \phantom{0} \notag \\
	\weight(\open \open\close \open L L \close \close) + \weight(\open \open L L \close \open\close \close) + \weight(\open \open L L \close \open L L \close \close).
\end{multline}

We need an analogue of \eqn{enlogicrule} for splitting a pattern $\open p_l p_r \close$ into the two subpatterns $p_l$ and $p_r$.  For this, examine each of the $2^l$ patterns that arose in partitioning $q$.  For each pattern $p = \open p_l p_r \close$ whose language is infinite (that is, the word $p$ contains the symbol $L$) and has $\depth(p) \geq d$, rewrite
\[
	\weight(p) = \weight(p_l) \cdot \weight(p_r) \cdot \mspace{-15mu} \prod_{\substack{0 \leq i \leq k \\ \text{$p$ matches $p_i$}}} \mspace{-20mu} x_{p_i},
\]
where `$p$ matches $p_i$' means that every tree in $L_p$ matches $p_i$ (so $L_p \subset L_{p_i}$).  If $L_p$ is infinite but $\depth(p) < d$, keep $\weight(p)$ intact as a variable.

Finally, for all tree patterns $p$ whose language is finite (i.e., $p$ is a tree), rewrite
\[
	\weight(p) = \prod_{0 \leq i \leq k} x_{p_i}^\text{number of copies of $p_i$ in $p$}.
\]

The algorithm is as follows.  As before, begin with the equation
\[
	\weight(L) = \weight(\open\close) + \weight(\open L L \close).
\]
At each step, take each new variable $\weight(q)$ and obtain another equation by performing the procedure described:  Write it as the sum of $2^l$ other variables, split the designated patterns into subpatterns, and explicitly compute the weights of any trees appearing.  Continue in this manner until there are no new variables produced; this must happen because we break up $\weight(p)$ whenever $\depth(p) \geq d$, so there are only finitely many possible variables.  Eliminate from this system of polynomial equations all but the $k+2$ variables $\weight(\open\close) = x_L$, $x_{p_1}$, \dots, $x_{p_k}$, and $\weight(L) = \En_{L, p_1, \dots, p_k}(x_L, x_{p_1}, \dots, x_{p_k})$ to obtain a polynomial equation satisfied by $\En_{L, p_1, \dots, p_k}(x_L, x_{p_1}, \dots, x_{p_k})$.

\section{Replacement bijections}\label{replacementbijections}

In this section we address the question of providing systematic bijective proofs of avoiding-equivalence.  Given two equivalent binary tree patterns $s$ and $t$, we would like to produce an explicit bijection between binary trees avoiding $s$ and binary trees avoiding $t$.  It turns out that this can often be achieved by structural replacements on trees.  We start by describing an example in full, and later generalize.

\subsection{An example replacement bijection}\label{examplereplacementbijection}

Consider the trees
\[
	t_2 = \vcentergraphics{binarytree4-2} \quad \text{and} \quad
	t_3 = \vcentergraphics{binarytree4-3}
\]
in class~4.2.  The idea is that since $n$-leaf trees avoiding $t_2$ are in bijection to $n$-leaf trees avoiding $t_3$, then somehow swapping all occurrences of these two tree patterns should produce a bijection.  However, since the patterns may overlap, it is necessary to specify an order in which to perform the replacements.  A natural order is to start with the root and work down the tree.  More precisely, a \emph{top-down replacement} is a restructuring of a tree $T$ in which we iteratively apply a set of transformation rules to subtrees of $T$, working downward from the root.

Take the replacement rule to be
\[
	\vcentergraphics{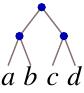} \to \vcentergraphics{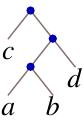},
\]
where the variables represent trees attached at the leaves, rearranged according to the permutation $3124$.  Begin at the root:  If $T$ itself matches the left side of the rule, then we restructure $T$ according to the rule; if not, we leave $T$ unchanged.  Then we repeat the procedure on the root's (new) children, then on their children, and so on, so that each vertex in the tree is taken to be the root of a subtree which is possibly transformed by the rule.  For example,
\[
	\vcentergraphics{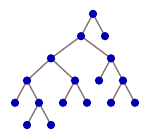} \to \vcentergraphics{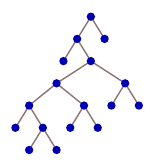} \to \vcentergraphics{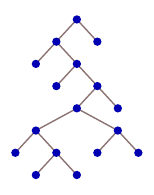} \to \vcentergraphics{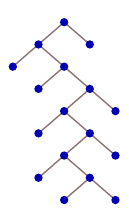}
\]
shows the three replacements required to compute the image (on the right) of a tree avoiding $t_2$.  The resulting tree avoids $t_3$.

This top-down replacement is invertible.  The inverse map is a \emph{bottom-up replacement} with the inverse replacement rule,
\[
	\vcentergraphics{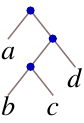} \to \vcentergraphics{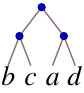}.
\]
Rather than starting at the root and working down the tree, we apply this map by starting at the leaves and working up the tree.

We now show that the top-down replacement is in fact a bijection from trees avoiding $t_2$ to trees avoiding $t_3$.  It turns out to be the same bijection given in \scn{4leafclasses} via words in $\{0,1\}^{n-2}$.

Assume $T$ avoids $t_2$; we show that the image of $T$ under the top-down replacement avoids $t_3$.  It is helpful to think of $T$ as broken up into (possibly overlapping) ``spheres of influence'' --- subtrees which are maximal with respect to the replacement rule in the sense that performing the top-down replacement on the subtree does not introduce instances of the relevant tree patterns containing vertices outside of the subtree.  It suffices to consider each sphere of influence separately.  A natural focal point for each sphere of influence is the highest occurrence of $t_3$.  We verify that restructuring this $t_3$ to $t_2$ under the top-down replacement produces no $t_3$ above, at, or below the root of the new $t_2$ in the image of $T$.

\begin{enumerate}
\item[above:] Since $t_3$ has depth $2$, $t_3$ can occur at most one level above the root of the new $t_2$ while overlapping it.  Thus it suffices to consider all subtrees with $t_3$ occurring at level $1$.  There are two cases,
\[
	\vcentergraphics{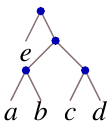} \quad \text{and} \quad
	\vcentergraphics{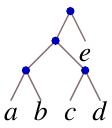}.
\]
The first case does not avoid $t_2$, so it does not occur in $T$.  The second case may occur in $T$.  However, we do not want the subtree itself to match $t_3$ (because we assume that the $t_3$ at level $1$ is the highest $t_3$ in this sphere of influence), so we must have $e = \open \close$.  Thus this subtree is transformed by the top-down replacement as
\[
	\vcentergraphics{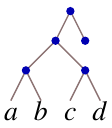} \to \vcentergraphics{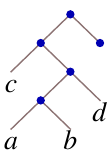}.
\]
The image does not match $t_3$ at the root, so $t_3$ does not appear above the root of the new $t_2$.
\item[at:] Since $T$ avoids $t_2$, every subtree in $T$ matching $t_3$ in fact matches the pattern $\open \open L L \close \open \open \close L \close \close$.  Such a subtree is restructured as
\[
	\vcentergraphics{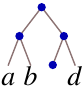} \to \vcentergraphics{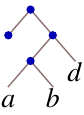}
\]
under the replacement rule, and the image does not match $t_3$ (because $c = \open \close$ is terminal).  Therefore the new $t_2$ cannot itself match $t_3$.
\item[below:] A general subtree matching $t_3$ and avoiding $t_2$ is transformed as
\[
	\vcentergraphics{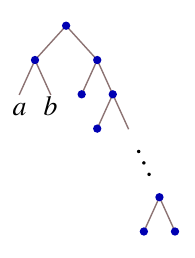} \to \vcentergraphics{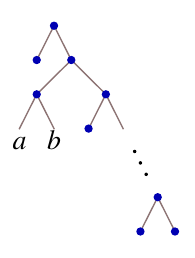} \to \vcentergraphics{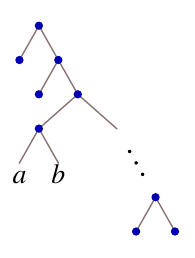} \to \; \cdots \; \to \vcentergraphics{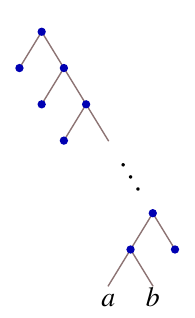} \to \; \cdots.
\]
Clearly $t_3$ can only occur in the image at or below the subtree $\open a b \close$.  However, since $\open a b \close$ is preserved by the replacement rule, any transformations on $\open a b \close$ can be considered independently.  That is, $\open a b \close$ is the top of a different sphere of influence, so we need not consider it here.  We conclude that $t_3$ does not occur below the root of the new $t_2$.
\end{enumerate}

If we already knew that $t_2$ and $t_3$ are equivalent (for example, by having computed $\Av_t(x)$ as in \scn{algorithm-avoidingsingle}), then we have now obtained a bijective proof of their equivalence.  Otherwise, it remains to show that if $T$ avoids $t_3$, then performing the bottom-up replacement produces a tree that avoids $t_2$; this can be accomplished similarly.

\subsection{General replacement bijections}\label{generalreplacementbijections}

A natural question is whether for any two equivalent binary tree patterns $s$ and $t$ there exists a sequence of replacement bijections and left--right reflections that establishes their equivalence.  For tree patterns of at most $7$ leaves the answer is ``Yes'', which perhaps suggests that these maps suffice in general.

\begin{conjecture}
If $s$ and $t$ are equivalent, then there is a sequence of top-down replacements, bottom-up replacements, and left--right reflections that produces a bijection from binary trees avoiding $s$ to binary trees avoiding $t$.
\end{conjecture}

In this section we discuss qualitative results regarding this conjecture.

Given two $m$-leaf tree patterns $s$ and $t$, one can ask which permutations of the leaves give rise to a top-down replacement that induces a bijection from trees avoiding $s$ to trees avoiding $t$.  Most permutations are not viable.  Candidate permutations can be found experimentally by simply testing all $m!$ permutations of leaves on a set of randomly chosen binary trees avoiding $s$; one checks that the image avoids $t$ and that composing the top-down replacement with the inverse bottom-up replacement produces the original tree.  This approach is feasible for small $m$, but it is slow and does not provide any insight into why certain trees are equivalent.  A question unresolved at present is to efficiently find all such bijections.

We \unabridged{now }return briefly to the replacement rule of \scn{examplereplacementbijection} to mention that a minor modification produces a bijection on the full set of binary trees.  Namely, take the two replacement rules
\[
	\vcentergraphics{binarytree4-3abcd} \to \vcentergraphics{binarytree4-2cabd} \quad \text{and} \quad
	\vcentergraphics{binarytree4-2abcd} \to \vcentergraphics{binarytree4-3bcad}.
\]
Again we perform a top-down replacement, now using both rules together.  That is, if a subtree matches the left side of either rule, we restructure it according to that rule.  Of course, it can happen that a particular subtree matches both replacement rules, resulting in potential ambiguity; in this case which do we apply?  Well, if both rules result in the same transformation, then it does not matter, and indeed with our present example this is the case.  To show this, it suffices to take the intersection $t_2 \cap t_3$ of the two left sides and label the leaves to represent additional branches that may be present:
\[
	\vcentergraphics{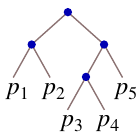}.
\]
Now we check that applying each of the two replacement rules to this tree produces the same labeled tree, namely
\[
	\vcentergraphics{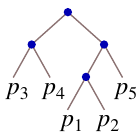}.
\]
Therefore we need not concern ourselves with which rule is applied to a given subtree that matches both.  Since the replacement rules agree on their intersection, the top-down replacement is again invertible and is therefore a bijection from the set of binary trees to itself.  By the examination of cases in \scn{examplereplacementbijection}, this bijection is an extension of the bijection between binary trees avoiding $t_2$ and binary trees avoiding $t_3$.

Thus we may choose from two types of bijection when searching for top-down replacement bijections that prove avoiding-equivalence.  The first type is from binary trees avoiding $s$ to binary trees avoiding $t$, using one rule for the top-down direction and the inverse for the bottom-up direction; these bijections in general do not extend to bijections on the full set of binary trees.  The second type is a bijection on the full set of binary trees, using both rules in each direction, that induces a bijection from binary trees avoiding $s$ to binary trees avoiding $t$.

We conclude with a curious example in which two tree patterns can only be proven equivalent by a two-rule bijection that does not involve them directly.  The trees
\[
	t_7 = \vcentergraphics{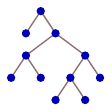} \quad \text{and} \quad
	t_{11} = \vcentergraphics{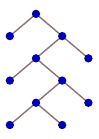}
\]
in class~6.5 are avoiding-equivalent by the permutation $126345$, but neither
\[
	t_{17} = \vcentergraphics{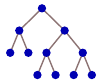}
\]
nor its left--right reflection has an equivalence-proving permutation to $t_7$, $t_{11}$, or their left--right reflections.  Thus, this equivalence cannot be established by a bijection that swaps $6$-leaf tree patterns.  However, it can be established by a bijection that swaps $4$-leaf tree patterns:  The previously mentioned bijection consisting of the two replacement rules
\[
	\vcentergraphics{binarytree4-3abcd} \to \vcentergraphics{binarytree4-2cabd} \quad \text{and} \quad
	\vcentergraphics{binarytree4-2abcd} \to \vcentergraphics{binarytree4-3bcad},
\]
induces a top-down replacement bijection from trees avoiding $t_7$ to trees avoiding $t_{17}$.  The reason is that $t_7$ and $t_{17}$ are formed by two overlapping copies of the class~4.2 trees
\[
	\vcentergraphics{binarytree4-2} \quad \text{and} \quad
	\vcentergraphics{binarytree4-3}
\]
respectively, and that $t_7$ and $t_{17}$ are mapped to each other under this bijection.

\section*{Acknowledgements}

I thank Phillipe Flajolet for helping me understand some existing literature, and I thank Lou Shapiro for suggestions which clarified some points in the paper.  Thanks to the referee for the reference to Stanley's book.

I am indebted to Elizabeth Kupin for much valuable feedback.  Her comments greatly improved the exposition and readability of the paper.  In addition, the idea of looking for one-rule bijections that do not extend to bijections on the full set of binary trees is hers, and this turned out to be an important generalization of the two-rule bijections I had been considering.

\appendix
\section*{Appendix. Table of equivalence classes}

This appendix lists equivalence classes of binary trees with at most $6$ leaves.  Left--right reflections are omitted for compactness.  For each class we provide a polynomial equation satisfied by $f = \En_{L, t}(x,y)$; an equation satisfied by $\Av_t(x)$ is obtained in each case by letting $y = 0$.

The data was computed by the \emph{Mathematica} package \textsc{TreePatterns} \cite{treepatterns} using Singular via the interface package by Manuel Kauers and Viktor Levandovskyy \cite{Singular interface package}.  Pre-computed data extended to $8$-leaf binary trees is now also available in \textsc{TreePatterns}.  The number of equivalence classes of $m$-leaf binary trees for $m = 1, 2, 3, \dots$ is $1, 1, 1, 2, 3, 7, 15, 44, \dots$.

\	

\subsubsection*{Class~1.1 \textnormal{(1 tree)}}
\[
	\includegraphics{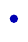}
\]
\footnotesize
\[
	x y f^2 - f + x y = 0
\]
\normalsize

\subsubsection*{Class~2.1 \textnormal{(1 tree)}}
\[
	\includegraphics{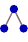}
\]
\footnotesize
\[
	x y f^2-f+x=0
\]
\normalsize

\subsubsection*{Class~3.1 \textnormal{(2 trees)}}
\[
	\includegraphics{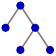}
\]
\footnotesize
\[
	x y f^2 + \left( - x^2 (y - 1) - 1\right) f + x = 0
\]
\normalsize

\subsubsection*{Class~4.1 \textnormal{(2 trees)}}
\[
	\includegraphics{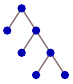}
\]
\footnotesize
\[
	\left(x y - x^3 (y - 1)\right) f^2 + \left( - x^2 (y - 1) - 1\right) f + x = 0
\]
\normalsize

\subsubsection*{Class~4.2 \textnormal{(3 trees)}}
\[
	\includegraphics{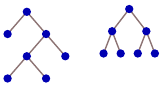}
\]
\footnotesize
\[
	x y f^2 + \left( - 2 x^2 (y - 1) - 1\right) f + \left(x^3 (y - 1) + x\right) = 0
\]
\normalsize

\subsubsection*{Class~5.1 \textnormal{(2 trees)}}
\[
	\includegraphics{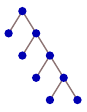}
\]
\footnotesize
\[
	- x^4 (y - 1) f^3 + \left(x y - x^3 (y - 1)\right) f^2 + \left( - x^2 (y - 1) - 1\right) f + x = 0
\]
\normalsize

\subsubsection*{Class~5.2 \textnormal{(10 trees)}}
\[
	\includegraphics{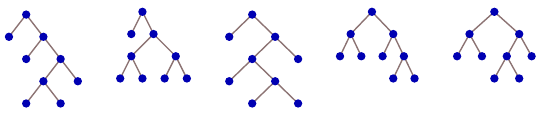}
\]
\footnotesize
\[
	\left(x y - x^3 (y - 1)\right) f^2 + \left(x^2 \left(x^2 - 2\right) (y - 1) - 1\right) f + \left(x^3 (y - 1) + x\right) = 0
\]
\normalsize

\subsubsection*{Class~5.3 \textnormal{(2 trees)}}
\[
	\includegraphics{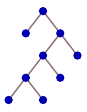}
\]
\footnotesize
\[
	x y f^3 + \left( - 3 x^2 (y - 1) - 1\right) f^2 + \left(3 x^3 (y - 1) + x\right) f - x^4 (y - 1) = 0
\]
\normalsize

\subsubsection*{Class~6.1 \textnormal{(2 trees)}}
\[
	\includegraphics{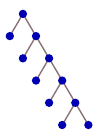}
\]
\footnotesize
\[
	- x^5 (y - 1) f^4 - x^4 (y - 1) f^3 + \left(x y - x^3 (y - 1)\right) f^2 + \left( - x^2 (y - 1) - 1\right) f + x = 0
\]
\normalsize

\subsubsection*{Class~6.2 \textnormal{(8 trees)}}
\[
	\includegraphics{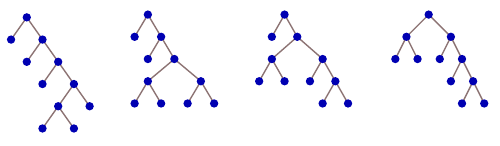}
\]
\footnotesize
\begin{multline*}
	- x^4 (y - 1) f^3 + x \left(x^2 \left(x^2 - 1\right) (y - 1) + y\right) f^2 + \phantom{0} \\
	\left(x^2 \left(x^2 - 2\right) (y - 1) - 1\right) f + \left(x^3 (y - 1) + x\right) = 0
\end{multline*}
\normalsize

\subsubsection*{Class~6.3 \textnormal{(14 trees)}}
\[
	\includegraphics{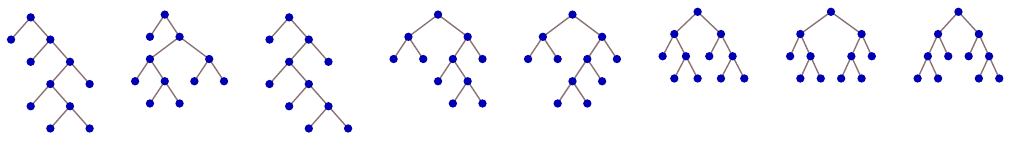}
\]
\footnotesize
\[
	x \left(x^2 \left(x^2 - 2\right) (y - 1) + y\right) f^2 + \left(2 x^2 \left(x^2 - 1\right) (y - 1) - 1\right) f + \left(x^3 (y - 1) + x\right) = 0
\]
\normalsize

\subsubsection*{Class~6.4 \textnormal{(8 trees)}}
\[
	\includegraphics{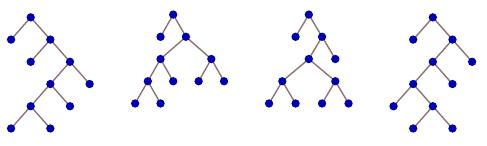}
\]
\footnotesize
\begin{multline*}
	\left(x y - x^3 (y - 1)\right) f^3 + \left(x^2 \left(2 x^2 - 3\right) (y - 1) - 1\right) f^2 + \phantom{0}\\
	\left( - x^5 (y - 1) + 3 x^3 (y - 1) + x\right) f - x^4 (y - 1) = 0
\end{multline*}
\normalsize

\subsubsection*{Class~6.5 \textnormal{(6 trees)}}
\[
	\includegraphics{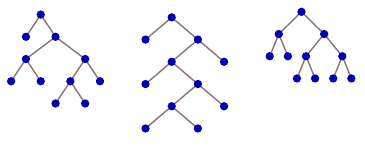}
\]
\footnotesize
\[
	\left(x y - 2 x^3 (y - 1)\right) f^2 + \left(x^2 \left(3 x^2 - 2\right) (y - 1) - 1\right) f + \left( - x^5 (y - 1) + x^3 (y - 1) + x\right) = 0
\]
\normalsize

\subsubsection*{Class~6.6 \textnormal{(2 trees)}}
\[
	\includegraphics{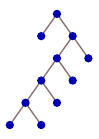}
\]
\footnotesize
\[
	- x y f^4 + \left(4 x^2 (y - 1) + 1\right) f^3 + \left( - 6 x^3 (y - 1) - x\right) f^2 + 4 x^4 (y - 1) f - x^5 (y - 1) = 0
\]
\normalsize

\subsubsection*{Class~6.7 \textnormal{(2 trees)}}
\[
	\includegraphics{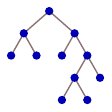}
\]
\footnotesize
\begin{multline*}
	x^4 \left(x^2 (y - 1) - y\right) (y - 1) f^3 + \left( - 2 x^7 (y - 1)^2 + x^5 (y - 1) (3 y - 2) - x^3 (y - 1) + x y\right) f^2 + \phantom{0} \\
	\left(x^2 \left(x^6 (y - 1) - 3 x^4 (y - 1) + x^2 - 2\right) (y - 1) - 1\right) f + \left(x^7 (y - 1)^2 + x^3 (y - 1) + x\right) = 0
\end{multline*}

\end{document}